\documentclass[11pt,a4paper]{article}
\pdfoutput=1 
\usepackage{graphicx}      

\usepackage{amsmath,amsfonts,amssymb,amscd}
\usepackage{bbm}
\usepackage[english]{babel}

\usepackage{graphicx}

\usepackage{tikz}

\def\IC{{\mathbb C}}

\def\IL{{\mathbb L}}

\newcommand{\sIL}{{{{\mathbb L}_s}}}
\newcommand{\csL}{{{{\cal L}_s}}}

\newcommand{\bA}{{\bf A}}
\newcommand{\bB}{{\bf B}}
\newcommand{\bC}{{\bf C}}

\newcommand{\bE}{{\bf E}}

\newcommand{\bL}{{\bf L}}
\newcommand{\bM}{{\bf M}}
\newcommand{\bN}{{\bf N}}
\newcommand{\bI}{{\bf I}}
\newcommand{\bH}{{\bf H}}

\newcommand{\bW}{{\bf W}}
\newcommand{\bR}{{\bf R}}
\newcommand{\bT}{{\bf T}}

\newcommand{\bx}{{\bf x}}

\newcommand{\bu}{{\bf u}}
\newcommand{\bV}{{\bf V}}

\newcommand{\bv}{{\bf v}}
\newcommand{\bw}{{\bf w}}

\newcommand{\bzeta}{{\bf \zeta}}
\newcommand{\bSigma}{{\bf \Sigma}}

\newcommand{\cO}{ {\cal O} }

\newcommand{\cL}{ {\cal L} }

\newcommand{\cR}{ {\cal R} }

\newcommand{\cV}{ {\cal V} }
\newcommand{\cW}{ {\cal W} }

\newcommand{\blambda}{\boldsymbol{\lambda}}
\newcommand{\bmu}{\boldsymbol{\mu}}

\newcommand{\bLambda}{\boldsymbol{\Lambda}}

\newtheorem{theorem}{Theorem}

\newtheorem{lemma}{Lemma}
\newtheorem{remark}{Remark}
\newtheorem{definition}{Definition}

\title{Loewner functions for bilinear systems} 

\author{Pauline Kergus\footnote{LAPLACE, Universit\'e de Toulouse, CNRS, INPT, UPS, Toulouse, France} \footnote{Corresponding author: pauline.kergus@cnrs.fr} , Ion Victor Gosea\footnote{Max Planck Institute for Dynamics of Complex Technical Systems, Sandtorstr. 1, 39106 Magdeburg, Germany} and Mihaly Petreczky\footnote{Centre de Recherche en Informatique, Signal et Automatique de Lille (CRIStAL), UMR CNRS 9189, CNRS Lille, France}}

\date{}
\begin{document}

\maketitle

\begin{abstract}                
This work brings together the moment matching approach based on Loewner functions and the classical Loewner framework based on the Loewner pencil in the case of bilinear systems. New Loewner functions are defined based on the bilinear Loewner framework, and a Loewner equivalent model is produced using these functions. This model is composed of infinite series that needs to be truncated in order to be implemented in practice. In this context, a new notion of approximate Loewner equivalence is introduced. In the end, it is shown that the moment matching procedure based on the proposed Loewner functions and the classical interpolatory bilinear Loewner framework both result in $\kappa$-Loewner equivalent models, the main difference being that the latter preserves bilinearity at the expense of a higher order.
\end{abstract}



\section{Introduction}

Accurate modeling of physical phenomena often leads to large-scale dynamical systems that require long simulation times and large data storage. In this context, model order reduction (MOR) aims at obtaining much smaller and simpler models that are still capable of accurately representing the behavior of the original process. Many different MOR techniques have been proposed: SVD-based (e.g., balanced truncation) methods, Krylov-based or moment-matching methods, proper orthogonal decomposition (POD) methods, and reduced basis methods. Most of these methods are included in the broad family of projection-based methods: the internal state $x$ is approximated by a vector of smaller dimension $\hat{x}$ obtained through a projection into a particular subspace. For more insights, the reader can refer to \cite{antoulas2000survey}, \cite{antoulas2005approximation},\cite{antoulas2010interpolatory}, \cite{baur2014model}, and to \cite{benner2015survey}.

Among them, the Loewner framework (LF) \cite{antoulas2017tutorial} is very appealing due to its data-driven nature, which makes it nonintrusive as it does not use the full/exact description of the model. It is based on the Loewner pencil that allows solving the generalized realization problem for linear time-invariant (LTI) systems \cite{mayo2007framework}, and obtaining reduced-order models through rational interpolation \cite{antoulas1986scalar}. In particular, the Loewner matrix is the cornerstone of the LF as it can be factored into the tangential generalized controllability matrix and the tangential generalized observability matrix, which can then be used to construct LTI models.

To broaden its applicability, extensions of the LF have been proposed for specific classes of nonlinear systems: bilinear \cite{morAntGI16}, quadratic-bilinear \cite{gosea2018data}, and linear switched systems \cite{gosea2018switched}. These extensions rely on higher-order transfer functions based on Volterra series, and the definition of associated observability and controllability matrices.

In parallel, an interconnection approach is proposed in \cite{simard2019interconnection}, based on the definition of left and right Loewner matrices. This approach allowed to develop a model order reduction procedure for linear time-varying systems \cite{simard2020loewner}. This idea is taken further in \cite{simard2021nonlinear} to extend the LF to general nonlinear input-affine systems by defining Loewner functions, which can be seen as time-varying nonlinear counterparts of the matrices classically used in the LF. This branch remains mostly theoretical since the explicit computation Loewner functions requires solving Partial Differential Equations (PDEs). 

This work aims at bringing these two interpretations of the Loewner framework together, in order to benefit from the practical implementations of \cite{morAntGI16} and the theoretical guarantees regarding steady-state behaviour coming with Loewner equivalence as in \cite{simard2021nonlinear}. To start with, this paper focuses on bilinear systems. New Loewner functions, inspired by \cite{morAntGI16}, are proposed for bilinear systems, allowing to derive a moment matching procedure for this category of systems. These Loewner functions are built based on the generalized reachability and observability matrices from \cite{morAntGI16}. In particular, they consist in infinite power series, which are not suitable for practical implementation. In this work, it is proposed to truncate these Loewner functions and a concept of approximate Loewner equivalence, namely \textit{$\kappa$-Loewner equivalence} is introduced. Finally, it is shown that the approach from \cite{morAntGI16} preserves bilinearity at the expense of a higher order than the moment matching approach resulting from the proposed Loewner function, though the latter results in a more complex model structure. Both models are $\kappa$-Loewner equivalent to the underlying bilinear full-order model, i.e., they both match the same first $\kappa$-th order moments.


The outline of this paper is as follows: in Section \ref{sec:preliminaries}, the interconnection-based approach of the LF for input-affine nonlinear systems as in \cite{simard2021nonlinear} and the LF for bilinear systems from \cite{morAntGI16} are recalled. Section \ref{sec:contribution} introduces Loewner functions for the case of bilinear systems, allowing to derive a reduced-order Loewner equivalent model. As the resulting model cannot be implemented in practice, it is approximated by truncating the Lowner functions, and the concept of \textit{$\kappa$-Loewner equivalence} is introduced. Connections between the Bilinear Loewner Framwork (BLF) from \cite{morAntGI16} and \cite{simard2021nonlinear} are then established.
Conclusions are provided in Section \ref{sec:conclusion}. 

\section{Preliminaries}
\label{sec:preliminaries}

\subsection{Nonlinear model order reduction based on Loewner functions}
\label{subsec:interconnectionLF}
\begin{figure*}[t]
    \centering
    \includegraphics[width=\textwidth]{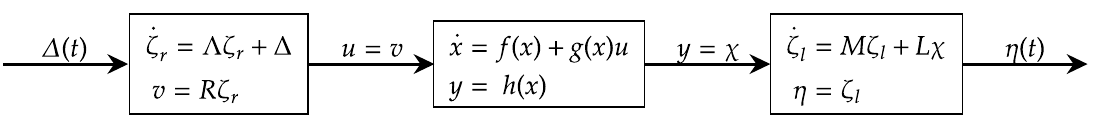}
    \caption{Interconnected system according to \protect\cite{simard2021nonlinear}.}
    \label{fig:interconnection}
\end{figure*}
In \cite{simard2021nonlinear}, nonlinear input-affine systems of the following form are considered:
\begin{equation}
    \begin{array}{r@{=}l}
        \Dot{\bx} & f(\bx) + g(\bx)u \\
        y & h(\bx)
    \end{array},
    \label{eq:nonlinear_system}
\end{equation}
with state $\bx\in \mathbb{C}^n$, input $u\in \mathbb{C}^m$ and output $y\in \mathbb{C}^p$.

This system is interconnected with two linear generators \eqref{eq:left_generator} and \eqref{eq:right_generator} as visible on Figure \ref{fig:interconnection}:
\begin{equation}
    \left\{\begin{array}{r@{=}l}
        \dot{\bzeta_r} & \bLambda\bzeta_r +\Delta \\
        v & \bR\bzeta_r
    \end{array}\right. ,
    \label{eq:left_generator}
\end{equation}
\begin{equation}
    \left\{\begin{array}{r@{=}l}
        \dot{\bzeta_l} & \bM\bzeta_l + \bL\chi \\
        \eta & \bzeta_l
    \end{array}\right. ,
    \label{eq:right_generator}
\end{equation}
with states $\bzeta_r\in\mathbb{C}^\rho$ and $\bzeta_l\in\mathbb{C}^\nu$, inputs $\Delta \in \mathbb{C}^\rho$ and $\chi \in \mathbb{C}^p$, and outputs $v\in\mathbb{C}^m$ and $v\in\mathbb{C}^\nu$, and with matrices $\bLambda\in\mathbb{C}^{\rho\times\rho}$, $\bR\in\mathbb{C}^{m\times\rho}$, $\bM\in\mathbb{C}^{\nu\times\nu}$, and $\bL\in\mathbb{C}^\nu\times p$. The matrices $\bLambda$ and $\bM$ have all eigenvalues on the imaginary axis, and these eigenvalues have geometric multiplicity one. For more details on the interconnection interpretation, see \cite{simard2019interconnection}.

Loewner functions are introduced as generalizations of the Loewner matrices to extend the Loewner method for model reduction to nonlinear input-affine systems. The tangential generalized controllability function $X : \mathbb{C}^\rho \rightarrow  \mathbb{C}^n$ is defined as the solution, provided it exists, to the PDE with boundary condition:
\begin{equation}
\begin{array}{lr}
    \frac{\partial X}{\partial \bzeta_r}\bLambda\bzeta_r = f(X(\bzeta_r))+g(X(\bzeta_r))\bR\bzeta_r  &    X(0) = 0
\end{array}
\label{eq:tang_gen_controllabillity}
\end{equation}
while the tangential generalized observability function $Y : \mathbb{C}^n \rightarrow  \mathbb{C}^\nu$ is defined as the solution, provided it exists, to the PDE with boundary condition
\begin{equation}
\begin{array}{lr}
    \frac{\partial Y}{\partial \bx}f(\bx) = \bM Y(\bx)-Lh(\bx)  &    Y(0) = 0
\end{array}.
\label{eq:tang_gen_observability}
\end{equation}
In \cite{simard2021nonlinear}, the existence of $X$ and $Y$ is proven under the following assumptions: the unforced system $\dot{\bx}=f(\bx)$ is locally exponentially stable at the origin and the matrices $\bLambda$ and $\bM$ have all eigenvalues on the imaginary axis and of geometric multiplicity one. The Loewner function $\cL$ is then defined as 
\begin{equation}
    \cL(\bzeta_r)=-Y(X(\bzeta_r)).
    \label{eq:Loewner_function}
\end{equation}
The nonlinear counterparts $\cV$ and $\cW$ of the data matrices $\bV$ and $\bW$ are defined as follows:
\begin{equation}
    \cV(\bzeta_r)=\frac{\partial Y}{\partial \bx}(X(\bzeta_r))g(X(\bzeta_r)), \cW(\bzeta_r)=h(X(\bzeta_r)).
    \label{eq:V_and_W}
\end{equation}
The left-Loewner function $\cL^l:\mathbb{C}^\rho\rightarrow\mathbb{C}^\nu$ is the solution (provided it exists) of
\begin{equation}
    \frac{\partial \cL^l}{\partial\bzeta_r}\bLambda\bzeta_r=\bM\cL^l(\bzeta_r)-\cV(\bzeta_r)\bR\bzeta_r, \cL^l(0)=0.
    \label{eq:left_loewner_function}
\end{equation}
The right-Loewner function  $\cL^r:\mathbb{C}^\rho\rightarrow\mathbb{C}^\nu$ is defined as:
\begin{equation}
    \cL^r(\bzeta_r)=\cL(\bzeta_r)-\cL^l(\bzeta_r).
    \label{eq:right_loewner_function}
\end{equation}
The shifted Loewner function $\csL:\mathbb{C}^\rho\rightarrow\mathbb{C}^\nu$ is defined as:
\begin{equation}
    \csL(\bzeta_r)=\bM\cL(\bzeta_r)+LW(\bzeta_r).
    \label{eq:shifted_loewner_function}
\end{equation}

\begin{definition}
    \label{def:loewner_equivalence}
    Two systems $\Sigma$ and $\overline{\Sigma}$ of the form \eqref{eq:nonlinear_system} are called \textit{Loewner equivalent} at $(\bLambda, \bR, \bM, \bL)$ if their left- and right-Loewner functions satisfy $\cL^l(\bzeta_r)=\overline{\cL}^l(\bzeta_r)$ and $\cL^r(\bzeta_r)=\overline{\cL}^r(\bzeta_r)$ in a neighborhood of the origin.
\end{definition}

The property of Loewner equivalence between two locally exponentially stable systems implies that, for initial conditions on the manifold $\bx=X(\zeta)$, the two systems interconnected with the generators have the same steady-state behavior, provided it exists.

\begin{theorem}
\label{theorem:loewner_eq_model}
    (From \cite{simard2021nonlinear}) Consider the interconnected system with $\rho=\nu$. Let $\cL^l$, $\cL^r$, $\cL$, $\csL$, $\cV$ and $\cW$ be the associated Loewner function. Assume that $\frac{\partial \cL}{\partial\bzeta_r}$ is non-singular. The system
    \begin{equation}
        \begin{array}{r@{=}l}
            \frac{\partial \cL}{\partial \bzeta_r}(\bx_r)\dot{\bx}_r & \csL(\bx_r)-\cV(\bx_r)u \\
            y_r & \cW(\bx_r)
        \end{array}
        \label{eq:loewner_eq_model}
    \end{equation}
    with state $\bx_r\in\mathbb{C}^\rho$, input $u\in\mathbb{C}^m$ and output $y_r\in\mathbb{C}^p$, is Loewner equivalent at $(\bLambda, \bR, \bM, \bL)$ to the system \eqref{eq:nonlinear_system}.
\end{theorem}

Locally, the original model and the interpolating model produce the same steady-state response, provided that it exists, when interconnected with generators corresponding to the Loewner functions. 

\subsection{The bilinear Loewner framework}
\label{subsec:bilinearLF}

The bilinear Loewner framework (BLF) presented here first appears in \cite{morAntGI16}. In this paper, single input - single output (SISO, $m=p=1$) bilinear systems of the form:
\begin{equation}
    \left\{\begin{array}{r@{=}l}
        \bE \dot{\bx} & \bA \bx + \bN \bx \bu + \bB u \\
        y & \bC \bx
    \end{array}\right.
    \label{eq:bilinear_system}
\end{equation} 
are considered, with state $\bx \in\mathbb{R}^n$, input $u \in \mathbb{R}$ and output $y \in \mathbb{R}$, with matrices $\bE,\bA,\bN\in\mathbb{R}^{n\times n}$, and $\bB, \bC^T\in\mathbb{R}^n$. 

The cornerstone of the BLF is that the behavior of bilinear systems is characterized by the following serie of generalized transfer function, defined for $l\geq 1$:
\begin{equation}
    \bH_l(s_1,s_2, \dots, s_l)=\bC\Phi(s_1)\bN\Phi(s_2)\bN\dots \bN\Phi(s_l)\bB,
    \label{eq:GTF}
\end{equation}
where 
\begin{equation}
    \Phi(s)=(s\bE-\bA)^{-1}
    \label{eq:resolvent}
\end{equation}
is the resolvent of the matrix pencil $(\bA,\bE)$. These generalized transfer functions are obtained by taking the Laplace transform of the kernels of the Volterra series expansion of the system, see \cite{morAntGI16} for more details.

The BLF interpolates the generalized transfer functions at interpolation points which are grouped into left and right multi-tuples, denoted $\bmu^{(j)}$ and $\blambda^{(i)}$ respectively, with $j=1 \dots \hat{q}$ and $i=1\dots \hat{k}$. The multi-tuples are defined as follows:
\small
\begin{equation}
    \bmu^{(j)}=\left\{\begin{array}{l}
        \left\{\mu_1^{(j)}\right\}  \\
        \left\{\mu_1^{(j)},\mu_2^{(j)}\right\}  \\
        \vdots \\
        \left\{\mu_1^{(j)},\mu_2^{(j)}, \dots , \mu_{p}^{(j)}\right\}  \\
    \end{array}\right. 
    \blambda^{(i)}=\left\{\begin{array}{l}
        \left\{\lambda_1^{(i)}\right\}  \\
        \left\{\lambda_2^{(i)},\lambda_1^{(i)}\right\}  \\
        \vdots \\
        \left\{\lambda_{m}^{(i)}, \dots ,\lambda_2^{(i)}, \lambda_1^{(i)}\right\}  \\
    \end{array}\right. .
    \label{eq:tuples}
\end{equation}
\normalsize


The generalized reachability matrix $\cR \in\mathbb{C}^{n\times k}$ 
associated with the right multi-tuples $\blambda^{(1)},\blambda^{(2)},\ldots,\blambda^{(\hat{k})}$ is
\begin{equation} \label{eq:reach_quad}
\cR=\left[\cR^{(1)},~\cR^{(2)},~\cdots,~\cR^{(\hat{k})}\right],
\end{equation}
where the matrices $\cR^{(j)} \in\mathbb{C}^{n\times m}$, $j = 1, \ldots, \hat{k}$, are associated with the $j$-th multi-tuple $\blambda^{(j)}$ in (\ref{eq:tuples}) are given by
\begin{equation}
    \begin{split}
        \cR^{(j)} = & \left[ \Phi(\lambda^{(j)}_1)\bB, \Phi(\lambda^{(j)}_2)\bN\Phi(\lambda^{(j)}_1)\bB, \dots \right. \\
        & \left. \dots , \Phi(\lambda^{(j)}_{m_i})\bN\dots \bN\Phi(\lambda^{(j)}_2)\bN\Phi(\lambda^{(j)}_1)\bB \right].
    \end{split}
    \label{eq:reach_Ri}
\end{equation}
Similarly, the generalized observability matrix  $\cO \in\mathbb{C}^{k \times n}$ associated with the left multi-tuples 
$\bmu^{(1)},\bmu^{(2)},\ldots,\bmu^{(\hat{k})}$ is given by
\begin{equation} \label{eq:obs_quad}
\cO=\left[  \big( \cO^{(1)} \big)^T,\; \big( \cO^{(2)} \big)^T, \; \ldots \big( \cO^{(\bar{k})} \big)^T \right]^T \in \IC^{k \times n},
\end{equation}
where $\cO^{(i)} \in\mathbb{C}^{p \times n}$, $i = 1, \ldots, \bar{k}$,
correspond to the $j$-th multi-tuple $\bmu^{(j)}$ in (\ref{eq:tuples}) and 
\begin{equation} \label{eq:obs_Oj}
\cO^{(j)}= \left[ \begin{array}{c} 
\bC^T\Phi(\mu^{(j)}_1) \\ 
\bC^T\Phi(\mu^{(j)}_1)\bN\Phi(\mu^{(j)}_2)\\
\vdots \\
\bC^T\Phi(\mu^{(j)}_1)\bN\Phi(\mu^{(j)}_2)\bN\dots\bN\Phi(\mu^{(j)}_{p_j})
\end{array} \right]. 
\end{equation}
Next, similarly to the linear case, the  Loewner matrix $\IL$ and the shifted Loewner matrix $\sIL$ are defined using the generalized
reachability (\ref{eq:reach_quad}) and observability (\ref{eq:obs_quad}) matrices as
\begin{equation} \label{eq:loewner_mat_quad}
\IL= - \cO\,\bE\,\cR,~~\sIL= - \cO\,\bA\,\cR \, .
\end{equation}
The fact that the Loewner matrices are factorized in terms of the pairs of matrices ($\bE,\bA$) and  $(\cO,\cR)$ is an inherent property of the Loewner framework, which holds true for both the bilinear and quadratic-bilinear case, and even for more general extensions involving higher-order polynomial structures.

Next, introduce the following matrices:
\begin{equation}
\bV = \cO \bB, \ \ \bW = \bC \cR, \ \ \bT = \cO \bN \cR,
\end{equation}

\begin{lemma}
    (from \cite{morAntGI16}) Assuming that the pencil $(\IL,\sIL)$ is regular, then the quintuple $$(-\IL,-\sIL,\bT,\bV,\bW)$$ defining the following model
    \begin{equation}
        \begin{array}{r@{=}l}
            -\IL \dot{\bx}_r & -\sIL \bx_r + \bT \bx_r u + \bV u \\
            y_r & \bW \bx_r
        \end{array}
        \label{eq:model_BLF0}
    \end{equation}
    is a minimal realization of an interpolant of the generalized transfer functions $\bH_l$ for $l=1\dots\hat{k}+\hat{q}$ at the points defined by the multi-tuples \eqref{eq:tuples}.
\end{lemma}

In case of redundant data, i.e. the pencil $(\IL,\sIL)$ is not regular, then a reduced-order model can be constructed by performing and truncating a singular value decomposition of the pencil, see \cite{morAntGI16}.
\section{Moment matching for bilinear systems}
\label{sec:contribution}
In this section, Loewner functions are proposed in Section \ref{subsec:loewner_fun} for bilinear systems as described in \eqref{eq:bilinear_system} with $\bE=I$. To start with, the left generator is given by $\bLambda =\lambda$ and $\bR=1$, and the right generator is given by $\bM=\mu$ and $\bL=1$, i.e. $\rho = \nu = 1$. Connections are then established with the BLF from \cite{morAntGI16} in Section \ref{subsec:connection_BLF}. The general case $\rho\geq1$ is presented in Section \ref{subsec:multiple_moments}. 

\subsection{Loewner functions for bilinear systems}
\label{subsec:loewner_fun}
\begin{theorem}\label{thm:rho1}
    The tangential generalized controllability function $X$
\begin{equation}
    X(\bzeta_r)=\sum_{k=1}^{\infty}\Phi_k \bzeta_r^k,
    \label{eq:controllability_fun}
\end{equation}
with coefficients $\left\{\Phi_k\right\}$ defined by
\begin{equation}
    \left\{\begin{array}{r@{=}l}
        \Phi_1 & (\lambda \bI -\bA)^{-1}\bB \\
        \Phi_k & (k\lambda \bI -\bA)^{-1}\bN\Phi_{k-1} \textnormal{ for } k>1 
    \end{array}\right.
    \label{eq:Phik}
\end{equation}
solves \eqref{eq:tang_gen_controllabillity}.
\end{theorem}

\textit{Proof:} First, $X(0)=0$. Note that, from \eqref{eq:Phik}, we have $k\lambda\Phi_k=\bA\Phi_k+\bN\Phi_{k-1}$ for $k>1$, and $\lambda\Phi_1=\bA\Phi_1+\bB$. Therefore, we have: 
$$\begin{array}{r@{=}l}
     \frac{\partial X}{\partial \bzeta_r} \lambda\bzeta_r & \sum_{k=1}^{\infty}k\lambda\Phi_k \bzeta_r^{k}\\
     & (\bA\Phi_1+\bB)\bzeta_r+ \sum_{k=2}^{\infty}(\bA\Phi_k + \bN\Phi_{k-1}) \bzeta_r^{k}\\
     & \bA\sum_{k=1}^{\infty}\Phi_k\bzeta_r^{k}+\left(\bN\sum_{k=1}^{\infty}\Phi_{k}\bzeta_r^{k}+\bB\right)\bzeta_r \\
     & f(X(\bzeta_r))+g(X(\bzeta_r))\bR\bzeta_r.
\end{array}$$

The tangential observability function
\begin{equation}
    Y(\bx)=C(\mu I-A)^{-1}\bx,
    \label{eq:observability_fun}
\end{equation}
solves \eqref{eq:tang_gen_observability} as $Y(0)=0$ and
$\frac{\partial Y}{\partial \bx } (\mu I-A)\bx = Cx$.


Then, the corresponding Loewner functions can be expressed as in \eqref{eq:Loewner_function}, \eqref{eq:V_and_W} and \eqref{eq:shifted_loewner_function} respectively:
\begin{equation}
    \cL(\bzeta_r)=-C(\mu I - A)^{-1}\sum_{k=1}^{\infty} \Phi_k \bzeta_r^k,
    \label{eq:loewner_function_bilinear}
\end{equation}
\begin{equation}
    \cV(\bzeta_r)=C(\mu I - A)^{-1}\left(B +\bN\sum_{k=1}^{\infty}\Phi_k\bzeta_r^k\right),
    \label{eq:V_bilinear}
\end{equation}
\begin{equation}
    \cW(\bzeta_r)=C\sum_{k=1}^{\infty}\Phi_k\bzeta_r^k,
    \label{eq:W_bilinear}
\end{equation}
\begin{equation}
    \csL(\bzeta_r)=- C(\mu I - A)^{-1}A\sum_{k=1}^{\infty} \Phi_k \bzeta_r^k.
    \label{eq:shifted_loewner_function_bilinear}
\end{equation}

In addition, the left and right Loewner functions are defined as:
\begin{equation}
    \cL^{l}(\bzeta_r)=C(\mu I - A)^{-1}\left(\frac{B}{\lambda-\mu}\bzeta_r+ \sum_{k=2}^{\infty} \frac{\bN\Phi_{k-1} }{(k \lambda - \mu)}\bzeta_r^k\right),
    \label{eq:left_Loewner_fun}
\end{equation}
\begin{equation}
    \begin{array}{rl}
        \cL^{r}(\bzeta_r) =& -C(\mu I - A)^{-1} \left (\Phi_1-\frac{B}{\lambda-\mu} \right) \bzeta_r  \\
         & -C(\mu I - A)^{-1}\sum_{k=2}^{\infty} (\Phi_k  - \frac{\bN \Phi_{k-1}}{(k \lambda - \mu)})\bzeta_r^k.
    \end{array}
    \label{eq:right_Loewner_fun}
\end{equation}

Then the Loewner equivalent model \eqref{eq:loewner_eq_model} at $(\lambda,1,\mu,1)$ to the system \eqref{eq:bilinear_system} is given by:
\begin{equation}
    \begin{array}{rl}
    \sum_{k=1}^{\infty} k C(\mu I-A)^{-1} \Phi_k r^{k-1} \dot{\bx}_r &  =\sum_{k=1}^{\infty} C(\mu I -A)^{-1}A\Phi_k r^k + \\
    & +\sum_{k=1}^{\infty} C(\mu I-A)^{-1}\bN\Phi_k r^k u_r  \\
    & + C(\mu I-A)^{-1} B u_r\\
        y_r =& \sum_{k=1}^{\infty} C\Phi_k r^k
    \end{array} 
   \label{eq:eq_model}
\end{equation}

\begin{remark}
    If $R\neq1$ and/or $L\neq 1$, the basis functions $\Phi_k$ are multiplied by $R$ in \eqref{eq:Phik} and the generalized observability function becomes $Y(\bx)=LC(\mu I-A)^{-1}\bx$.
\end{remark}

\subsection{Truncation of the Loewner functions and $\kappa$-Loewner equivalence}
\label{subsec:truncation}
In practice, it is difficult to implement the Loewner equivalent model \eqref{eq:eq_model} as the Loewner functions consist in infinite series. To overcome this issue, let the truncated (generalized) controllability function be defined for a given $\kappa$ as:
\begin{equation}
\hat{X}(\bzeta_r)=\sum_{k=1}^{\kappa}\Phi_k \bzeta_r^k.
\label{eq:controllability_fun_trunc}
\end{equation}
As a result, the truncated left and right Loewner functions are given by
\begin{equation}
    \hat{\cL}^{l}(\bzeta_r)=C(\mu I - A)^{-1}\left(\frac{B}{\lambda-\mu}\bzeta_r+ \sum_{k=2}^{\kappa+1} \frac{\bN\Phi_{k-1} }{(k \lambda - \mu)}\bzeta_r^k\right),
    \label{eq:left_Loewner_fun_trunc}
\end{equation}
and
\begin{equation}
    \begin{array}{rl}
        \hat{\cL}^{r}(\bzeta_r) = & -C(\mu I - A)^{-1} \left (\Phi_1-\frac{B}{\lambda-\mu} \right) \bzeta_r  \\
         & -C(\mu I - A)^{-1}\sum_{k=2}^{\kappa} (\Phi_k  - \frac{\bN \Phi_{k-1}}{(k \lambda - \mu)})\bzeta_r^k \\
         & -C(\mu I - A)^{-1}\frac{\bN \Phi_{\kappa}}{((\kappa+1) \lambda - \mu)})\bzeta_r^{\kappa+1} .
    \end{array}
    \label{eq:right_Loewner_fun_trunc}
\end{equation}

Hence the truncated model \eqref{eq:loewner_eq_model} is given by:
\begin{equation}
    \begin{array}{rl}
    \sum_{k=1}^{\kappa} k C(\mu I-A)^{-1} \Phi_k r^{k-1} \dot{\bx}_r &  =\sum_{k=1}^{\kappa} C(\mu I -A)^{-1}A\Phi_k r^k + \\
    & +\sum_{k=1}^{\kappa} C(\mu I-A)^{-1}\bN\Phi_k r^k u_r  \\
    & + C(\mu I-A)^{-1} B u_r\\
        y_r =& \sum_{k=1}^{\kappa} C\Phi_k r^k
    \end{array} 
   \label{eq:eq_model_truncated}
\end{equation}

The truncation introduced previously for practical implementation is supported by the following concept of $\kappa$-Loewner equivalence. Intuitively, $\kappa$-Loewner equivalence means that the Loewner functions are close for sufficiently small signals $\zeta_r$.

\begin{definition}
\label{kappa_Loewner_eq}
(\textbf{$\kappa$-Loewner equivalence})
    Two systems $\bSigma$ and $\hat{\bSigma}$ are called $\kappa$-Loewner equivalent at $(\bLambda,L,\bM, R)$, if the left- and right- Loewner functions of $\bSigma$ and $\hat{\bSigma}$, denoted by $\mathcal{L}^{r}(\zeta_r)$,
$\mathcal{L}^{l}(\zeta_r)$ and $\hat{\mathcal{L}}^{r}(\zeta_r)$, $\hat{\mathcal{L}}^{l}(\zeta_r)$ respectively, are smooth and have the same $k$-order derivatives at zero for any $k \leq \kappa$:
\begin{equation}
\left\{
    \begin{array}{l}
       \frac{d^k}{d\zeta_r^k} \mathcal{L}^{r}(\zeta_r)|_{\zeta_r=0}=\frac{d^k}{d\zeta_r^k}\hat{\mathcal{L}}^{r}(\zeta_r)|_{\zeta_r=0} \\
        \frac{d^k}{d\zeta_r^k} \mathcal{L}^{l}(\zeta_r)|_{\zeta_r=0}=\frac{d^k}{d\zeta^k_r} \hat{\mathcal{L}}^{l}(\zeta_r)|_{\zeta_r=0} \\
    \end{array} 
    \right.
    \label{eq:kappa_loewner_eq}
\end{equation}
\end{definition}

Note that Loewner equivalence corresponds to the case $\kappa=\infty$.

\begin{theorem}
    Given $\lambda$, $\mu$, the truncated model \eqref{eq:eq_model_truncated} is $\kappa$-Loewner equivalent to the original system at $(\lambda,1,\mu, 1)$.
    \label{thm:kappa_Loewner_model}
\end{theorem}
\textit{Proof:} Considering the left and right Loewner functions in \eqref{eq:left_Loewner_fun} and \eqref{eq:right_Loewner_fun}, their truncated counterparts in \eqref{eq:left_Loewner_fun_trunc} and \eqref{eq:right_Loewner_fun_trunc} satisfy \eqref{eq:kappa_loewner_eq}.

\subsection{Connection to the bilinear Loewner framework}
\label{subsec:connection_BLF}
In order to establish a first connection between the moment matching method by \cite{simard2021nonlinear} for bilinear systems with the BLF from \cite{morAntGI16}, the following interpolation multi-tuples are considered:
\begin{equation}
    \bmu^{(1)}=\left\{\mu\right\} \ \text{and} \ 
    \blambda^{(1)}=\left\{\begin{array}{l}
        \left\{\lambda\right\}  \\
        \left\{2\lambda,\lambda \right\}  \\
        \vdots \\
        \left\{\kappa\lambda, \dots ,2\lambda, \lambda \right\}  \\
    \end{array}\right. .
    \label{eq:tuples_moment_matching}
\end{equation}
The reachability matrix $\cR$ from the BLF can be written as follows: \begin{equation}
    \cR = \begin{bmatrix}
       \Phi_1 & \Phi_2 & \cdots & \Phi_{\kappa} 
    \end{bmatrix} \in \mathbb{C}^{n\times \kappa },
    \label{eq:R_BLF}
\end{equation}
and the observability matrix is $\cO = \bC(\mu\bI-\bA)^{-1}\in\mathbb{C}^{1\times n}$.

Then, the generalized controllability function can be written as
\begin{equation}
    X(\bzeta_r)= \mathcal{\bR}\bv(\bzeta_r),
    \label{eq:linkBLF_cont}
\end{equation}
with $\bv(\bzeta_r)\in\mathbb{R}^{\kappa}$ defined as:
\begin{equation}
    \bv(\bzeta_r)= \left[\zeta_{r},\zeta_{r}^2, \dots, \zeta_{r}^\kappa\right].
    \label{eq:v_zeta}
\end{equation}
The generalized observability function can be written as
\begin{equation}
    Y(\bx)= \mathcal{O}\bx.
    \label{eq:linkBLF_obs}
\end{equation}
The Loewner functions can then be expressed according to the Loewner matrices as follows: 
\begin{equation}
\begin{split}
& \cL(\bzeta_r)=-\cO\cR\bv(\bzeta_r)=-\IL\bv(\bzeta_r),\\
& \cV(\bzeta_r)=\cO\left(B +\bN\cR\bv(\bzeta_r)\right)=\bV+\bT\bv(\bzeta_r),\\
& \cW(\bzeta_r)=C\cR\bv(\bzeta_r)=\bW\bv(\bzeta_r),\\
& \csL(\bzeta_r)=-\bM\cO\cR\bv(\bzeta_r)+LC\cR\bv(\bzeta_r)=\sIL\bv(\bzeta_r).
\end{split}
\label{eq:loewner_functions_BLF}
\end{equation}
The BLF then results in the following bilinear model of order $\kappa$:
\begin{equation}
    \begin{array}{r@{=}l}
        -\IL \dot{\bx}_r & -\sIL \bx_r - \bT \bx_r u - \bV u \\
        y_r & \bW x_r
    \end{array}
    \label{eq:model_BLF}
\end{equation}
while the Loewner equivalent model obtained based on the Loewner functions is of order 1 and is characterized by:
\begin{equation}
    \begin{array}{r@{=}l}
        -\IL\frac{\partial\bv}{\partial \bzeta_r}(\bx_r)\dot{\bx}_r & \sIL\bv(\bx_r)-(\bT\bv(\bx_r)+\bV)u \\
        y_r & \bW\bv(\bx_r)
    \end{array}
    \label{eq:loewner_eq_model_rho1}
\end{equation}

Note that when $\kappa=1$, the Loewner equivalent model \eqref{eq:loewner_eq_model_rho1} 
is the same than the one resulting from the BLF \eqref{eq:model_BLF}, and is written in descriptor form.

For $\kappa\geq 2$, the BLF results in a bilinear model of order $\kappa$ while the moment matching approach results in a first-order model but that does not preserve bilinearity. 

\begin{remark}
    Note that for $\kappa\geq2$, the model \eqref{eq:loewner_eq_model_rho1} from the BLF is not a descriptor model, but an ordinary differential equation of order $\kappa$. This is due to the choice of the interpolation multi-tuples. In the classical BLF, the data is usually split equally between left and right multi-tuples, and the left and right multi-tuples have the same number of components, which leads to descriptor models (with square Loewner matrices).
\end{remark}




\begin{theorem}
    The BLF model \eqref{eq:model_BLF} is $\kappa$-Loewner equivalent to the original one \eqref{eq:bilinear_system}.
\end{theorem}
\textit{Proof:} The Loewner functions corresponding to the BLF models are derived as in the moment matching procedure, based on the generalized controllability and observability functions which can be written as follows:
\begin{equation}
    \begin{array}{l}
        Y_{BLF}(x)=\bW(-\mu\IL+\sIL)x , \\
        X_{BLF}(\bzeta_r)=\sum_{k=1}^{\kappa}\Phi_k^{BLF}\bzeta_r^k,
    \end{array}
    \label{eq:Loewner_fct_BLF}
\end{equation}
with
$$\left\{\begin{array}{r@{=}l}
\Phi^{BLF}_1 & (-\lambda \IL +\sIL)^{-1}\bV\\
\Phi_k^{BLF} & (-k\lambda_i \IL +\sIL)^{-1}\bT\Phi_{k-1}^{BLF} \textnormal{ for } k>1 
\end{array}\right. .$$
The rest of the Loewner functions can be derived accordingly, in a similar manner as to \eqref{eq:loewner_function_bilinear}-\eqref{eq:right_Loewner_fun}.

Due to the interpolatory properties of the BLF, see \cite{morAntGI16}, we have that, for all $l=1\dots \kappa$:
\begin{equation}
    \begin{array}{l}
        \bH_l(\lambda,2\lambda, \dots, l\lambda)= \bH_l^{BLF}(\lambda,2\lambda, \dots, l\lambda) \\
        \bH_{l+1}(\mu,\lambda,2\lambda, \dots, l\lambda)=\bH_{l+1}^{BLF}(\mu,\lambda,2\lambda, \dots, l\lambda) \\
    \end{array}
    \label{eq:BLF_interpolation}
\end{equation}
where $\bH_l^{BLF}$ are the generalized transfer functions of the BLF model, as defined in \eqref{eq:GTF}. According to the definition of the generalized transfer functions in \eqref{eq:GTF} and to the expressions of the Loewner functions in the bilinear case, it implies that the truncated Loewner functions corresponding to the BLF model \eqref{eq:model_BLF} are equal to the truncated ones from the moment matching procedure proposed in Section \ref{subsec:loewner_fun}. Consequently, according to Theorem \ref{thm:kappa_Loewner_model}, the BLF model is $\kappa$-Loewner equivalent to the original system.

\subsection{Extension to higher-order generators}
\label{subsec:multiple_moments}
These Loewner functions can be extended to the case $\rho\geq1$ (recall that $\rho=\nu$). Consider the general case where the left generator is given by $\bLambda$ and $\bR$, see \eqref{eq:left_generator}, and the right generator is given by $\bM$ and $\bL$, see \eqref{eq:right_generator}.  Up to a similarity transformation, the left generator can be defined by a diagonal matrix $\bLambda$, with distinct coefficients $\lambda_1, \dots,\lambda_\rho$ on the imaginary axis and $\bR=\left[R_1 \ \dots \ R_\rho\right]$. The same goes for the right generator, and the coefficients of $\bM$ are denoted $\mu_1, \dots,\mu_\rho$ and $\bL^T=\left[L_1 \ \dots \ L_\rho\right]$.

\begin{theorem}
    For $\bLambda\in\mathbb{C}^{\rho\times\rho}$, $\bR^T\in\mathbb{C}^\rho$, with $\rho\geq 1$, the tangential generalized controllability function $X$ is defined by
\begin{equation}
    X(\bzeta_r)=\left[X^{(1)}(\zeta_{r_1})\ \dots \ X^{(\rho)}(\zeta_{r_\rho})\right],
    \label{eq:controllability_fun_gen}
\end{equation}
where, for $i=1\dots\rho$, $\zeta_{r_i}$ denotes the $i$-th component of $\bzeta_r$, and the functions $X^{(i)}$ are defined as in the first-order generator case as follows:
\begin{equation}
    X^{(i)}(\bzeta_r^{(i)})=\sum_{k=1}^{\infty}\Phi_k^{(i)} \zeta_{r_i}^k
    \label{eq:Xi}
\end{equation}
with the coefficients $\left\{\Phi_k^{(i)}\right\}$ defined by
\begin{equation}
    \left\{\begin{array}{r@{=}l}
        \Phi_1^{(i)} & (\lambda_i \bI -\bA)^{-1}\bB R_i\\
        \Phi_k^{(i)} & (k\lambda_i \bI -\bA)^{-1}\bN\Phi_{k-1}^{(i)}R_i \textnormal{ for } k>1 
    \end{array}\right. .
    \label{eq:Phiki}
\end{equation}
The tangential generalized controllability function $X$ from \eqref{eq:controllability_fun_gen} solves \eqref{eq:tang_gen_controllabillity}.
\end{theorem}

\textit{Proof:} First, as $X^{(i)}(0)=0$ for all $i=1\dots \rho$, we have $X(0)=0$. In addition, from Theorem \ref{thm:rho1}, we have, for all $i=1\dots \rho$, that:
$$\frac{\partial X^{(i)}}{\partial \zeta_{r,i}} \lambda_i\zeta_{r,i} = f(X^{(i)}(\zeta_{r,i}))+g(X^{(i)}(\zeta_{r,i}))R_i\zeta_{r,i}.$$
which implies that $X$ from \eqref{eq:controllability_fun_gen} solves \eqref{eq:tang_gen_controllabillity}.

The generalized observability function can be defined in a similar way.
\begin{theorem}
The tangential observability function is defined as:
\begin{equation}
    Y(\bx)=\begin{bmatrix}
        Y^{(1)}(\bx) \\ \vdots \\ Y^{(\rho)}(\bx)
    \end{bmatrix},
    \label{eq:observability_fun_gen}
\end{equation}
with, for $i=1\dots\rho$, the function $Y^{(i)}$ defined as in the first-order generator case as $Y^{(i)}(\bx)=L_i\bC(\mu_iI-\bA)^{-1}\bx$. The tangential generalized controllability function $X$ from \eqref{eq:observability_fun_gen} solves \eqref{eq:tang_gen_observability}.
\end{theorem}

\textit{Proof:} First, as $Y^{(i)}(0)=0$ for all $i=1\dots \rho$, we have $Y(0)=0$. In addition, from Section \ref{subsec:loewner_fun}, we have for all $i=1\dots \rho$:
$$\frac{\partial Y^{(i)}}{\partial \bx } f(\bx) = \mu_iY^{(i)}\bx-L_ih(\bx)$$
which implies that $Y$ from \eqref{eq:observability_fun_gen} solves \eqref{eq:tang_gen_observability}.

Similarly to the first-order generators case, a truncated version of the tangential generalized controllability function can be introduced so that the resulting model can be implemented in practice:
\begin{equation}
    \forall i=1\cdots\rho, \ \hspace{0.2cm} \ \hat{X}^{(i)}(\bzeta_r^{(i)})=\sum_{k=1}^{\kappa}\Phi_k^{(i)} \zeta_{r_i}^k.
    \label{eq:Xi_truncated}
\end{equation}
For any value of $\kappa$ (finite or infinite), the connection with the BLF from \cite{morAntGI16} can be established as in Section \ref{subsec:connection_BLF} by choosing $\rho$
interpolation multi-tuples:
\begin{equation}
    \bmu^{(i)}=\left\{\mu_i\right\} \ \text{and} \ 
    \blambda^{(i)}=\left\{\begin{array}{l}
        \left\{\lambda_i\right\}  \\
        \left\{2\lambda_i,\lambda_i \right\}  \\
        \vdots \\
        \left\{\kappa\lambda_i, \dots ,2\lambda_i, \lambda_i \right\}  \\
    \end{array}\right. .
    \label{eq:tuples_moment_matching_gen}
\end{equation}
Then, equation \eqref{eq:linkBLF_cont} becomes:
\begin{equation}
\begin{split}
& \cL(\bzeta_r)=-\IL\bw(\bzeta_r),\\
& \cV(\bzeta_r)=\bV+\bT\bw(\bzeta_r),\\
& \cW(\bzeta_r)=\bW\bw(\bzeta_r),\\
& \csL(\bzeta_r)=-\sIL\bw(\bzeta_r),
\end{split}
\label{eq:loewner_functions_BLF_gen}
\end{equation}
with $\bw(\bzeta_r)\in\mathbb{R}^{\rho\kappa\times \rho}$ defined as:
\begin{equation*}
    \bw(\bzeta_r)= 
    \begin{bmatrix}
        \bv(\zeta_{r,1}) & 0 & \cdots & 0\\
        0 & \bv(\zeta_{r,2}) &  & \vdots \\
        \vdots &  & \ddots & 0 \\
        0 & \cdots & 0 & \bv(\zeta_{r,\rho})\\
    \end{bmatrix}
\end{equation*}
with $\bv$ defined in \eqref{eq:v_zeta}.

In the end, the BLF results in a bilinear model of order $\rho\kappa$ while the moment matching approach results in a model of order $\rho$ that does not preserve bilinearity. As previously shown, the two models are exactly the same for $\kappa=1$. Finally, the two models are $\kappa$-Loewner equivalent to the original one at $(\bLambda,\bL,\bM,\bR)$.

\section{Conclusion}
\label{sec:conclusion}
In this work, new Loewner functions have been proposed for bilinear systems, using the generalized reachability and observability matrices used in the bilinear Loewner framework from \cite{morAntGI16}. A reduced-order Loewner equivalent model is then derived as in \cite{simard2021nonlinear}. Inpractice, it is quite challenging to implement. That is why, in this work, the Loewner functions are truncated in order to obtain approximate Loewner equivalence.
This is defined in this work as \textit{$\kappa$-Loewner equivalence}. It appears that the approach from \cite{morAntGI16} preserves bilinearity at the expense of a higher order than the moment matching approach resulting from the proposed Loewner function, although the latter results in a more complex model structure. Regarding the choice of $\kappa$, practical implementations of the BLF from \cite{morAntGI16} are performed with $\kappa=2$ but based on different multi-tuples definition. Future work will investigate the impact of $\kappa$ and could provide approximation error bounds for the truncated $\kappa$-Loewner equivalent model. In addition, the authors aim at expanding this concept and proposing Loewner functions in a similar manner for more general polynomial extensions of the Loewner framework.




\bibliographystyle{unsrt}             
\bibliography{ifacconf}

\end{document}